\documentclass[12pt]{amsart}
\usepackage{amsmath}
\usepackage{amsthm}
\usepackage{amssymb}
\usepackage{eucal}
\usepackage{times}
\usepackage{euler}
\usepackage{eufrak}

\author{Piotr \'Sniady}
\thanks{Research supported by State Committee for Scientific
Research (KBN) grant \mbox{2 P03A 007 23}.}

\address{Institute of Mathematics,
University of Wroclaw, pl.~Grunwaldzki~2/4, 50-384 Wroclaw,
Poland}
\email{Piotr.Sniady@math.uni.wroc.pl}
\title{Factoriality of Bo\.zejko--Speicher von Neumann algebras}

\theoremstyle{plain}
\newtheorem{lemma}{Lemma}[section]
\newtheorem{theorem}[lemma]{Theorem}
\newtheorem{proposition}[lemma]{Proposition}

\theoremstyle{definition}

\theoremstyle{remark}
\newtheorem*{remark}{Remark}

\newcommand{\Ha}{\mathcal{H}}
\newcommand{\Hado}[1]{\Ha^{\otimes #1}}

\newcommand{\C}{{\mathbb{C}}}
\newcommand{\R}{{\mathbb{R}}}

\newcommand{\F}{{\mathcal{F}}}
\newcommand{\gwia}{^{\ast}}

\newcommand{\cc}{M}
\newcommand{\can}{m}
\newcommand{\ccre}{m^{\dag}}

\newcommand{\leftqan}{l}
\newcommand{\rightqan}{r}
\newcommand{\leftfreecre}{l^{\ast}}
\newcommand{\rightfreecre}{r^{\ast}}
\newcommand{\leftqcre}{l^{\ast}}
\newcommand{\rightqcre}{r^{\ast}}
\newcommand{\leftfield}{L}
\newcommand{\rightfield}{R}

 \DeclareMathOperator{\inv}{inv}
 
 \DeclareMathOperator{\dimm}{dim}

\begin{document}

\maketitle

\begin{abstract}
We study the von Neumann algebra generated by $q$--deformed
Gaussian elements $\leftqan_i+\leftqcre_i$ where operators
$\leftqan_i$ fulfill the $q$--deformed canonical commutation
relations $\leftqan_i \leftqcre_j-q \leftqcre_j
\leftqan_i=\delta_{ij}$ for $-1<q<1$. We show that if the number
of generators is finite, greater than some constant depending on
$q$, it is a $\text{II}_1$ factor which does not have the property
$\Gamma$. Our technique can be used for proving factoriality of
many examples of von Neumann algebras arising from some
generalized Brownian motions, both for type $\text{II}_1$ and type
$\text{III}$ case.
\end{abstract}

\section{Introduction}

In 1970 Frisch and Bourret \cite{FrischBourret} considered
commutation relation
\begin{equation}
\label{eq:qdeform}
\leftqan_\phi \leftqcre_\psi - q \
\leftfreecre_\psi \leftqan_\phi = (\phi,\psi)
\end{equation}
fulfilled by annihilation operator $\leftqan$ and its adjoint
creation operator $\leftqcre$ for $q\in\R$. Annihilation and
creation operators are indexed by elements of a fixed real Hilbert
space $\Ha_\R$ and they act on a complex Hilbert space $\F$,
called Fock space. There is a distinguished unital vector
$\Omega\in\F$, called vacuum, such that $\leftqan_\phi \Omega=0$
for every $\phi\in\Ha_\R$. It was a long time until Bo\.zejko and
Speicher \cite{BozejkoSpeicher1994} showed in 1994 the existence
of the operators considered by Frisch and Bourret for all $-1< q
<1$.

Frisch and Bourret were studying generalized Gaussian variables
$L_\phi=\leftqan_\phi+\leftqcre_\phi$; it turns out that the
algebra $\Gamma_q$ generated by $(L_\phi)$ can be equipped with a
tracial state $a\mapsto \langle \Omega, a \Omega\rangle$. The
motivation for studying such Gaussian variables is that if $q=1$
then \eqref{eq:qdeform} coincides with the canonical commutation
relations and hence $(L_\phi)$ can be identified with a family of
classical Gaussian random variables; for $q=-1$ relation
\eqref{eq:qdeform} coincides with the canonical anticommutation
relations; furthermore it turned out much later that for $q=0$
relation \eqref{eq:qdeform} coincides with the free relation and
hence $(L_\phi)$ is a family of Voiculescu semicircular elements
\cite{VoiculescuDykemaNica}.

In this article we will study the von Neumann algebra $\Gamma_q$
generated by $q$--deformed Gaussian variables $L_\phi$. Since for
$q=0$ these von Neumann algebras are isomorphic to the free group
factors, we can consider the general case as a `smooth'
deformation of this eminent case. Not too much is known about
$\Gamma_q$, in particular it is not clear if it is always
isomorphic to the free group factors. Bo\.zejko and Speicher
\cite{BozejkoSpeicher1994} showed that under certain conditions
$\Gamma_q$ is non--injective; recently Nou \cite{Nou} showed that
it is enough to assume that $-1<q<1$ and $\dimm \Ha_\R\geq 2$.
Recently Shlyakhtenko \cite{Shlyakhtenko2003} showed that if $|q|<
\sqrt{2}-1$ and $\dimm \Ha_\R\geq 2$ then algebras $\Gamma_q$ are
solid (cf.\ \cite{OzawaSolid}) and if they are factors then they
do not have the property $\Gamma$.

Bo\.zejko, K\"ummerer and Speicher \cite{BozejkoKummererSpeicher}
showed that if the number of generators is infinite (i.e.\ if
$\dimm \Ha_\R=\infty$) then for every $-1<q<1$ the algebra
$\Gamma_q$ is a $\text{II}_1$ factor. In this article we show that
this result remains true if the number of generators is finite,
greater than some constant depending on $q$ and that this factor
does not have the property $\Gamma$. We also point out that the
same proof can be used for proving the analogous result for many
other von Neumann algebras, both finite and infinite.

\section{Notations}
If not stated otherwise, all results presented in this section are
due to Bo\.zejko and Speicher \cite{BozejkoSpeicher1994}.

\subsection{Fock space}
Let $\Ha_\R$ be a real Hilbert space equipped with a bilinear
scalar product $(\cdot,\cdot)$; we will denote its
complexification by $\Ha$ and the corresponding sesquilinear
scalar product by $\langle\cdot,\cdot\rangle$. Let furthermore
$-1<q<1$ be fixed. For integer $n\geq 0$ we introduce an operator
$P^{(n)}:\Hado{n} \rightarrow \Hado{n}$ given by
\begin{equation}
\label{eq:symmetrizer} P^{(n)}(\psi_1\otimes\cdots\otimes
\psi_n)=\sum_{\sigma\in S_n} q^{\inv\sigma}\
\psi_{\sigma(1)}\otimes\cdots\otimes \psi_{\sigma(n)},
\end{equation} where $\inv\sigma$ is the number of inversions in
$\sigma$, i.e.~the number of pairs $(i,j)$ such that $1\leq
i<j\leq n$ and $\sigma(i)>\sigma(j)$. Operator $P^{(n)}$ is
strictly positive, therefore we can equip $\Ha^{\otimes n}$ with a
new scalar product
$$ \langle \Phi, \Psi \rangle_q = \langle \Phi, P^{(n)} \Psi
\rangle,$$ for $\Phi,\Psi\in\Ha^{\otimes n}$ where $\langle
\cdot,\cdot \rangle$ denotes the standard scalar product on
$\Ha^{\otimes n}$. In the following by $\Ha^{\otimes n}$ we will
mean the $n$--fold tensor product equipped with the standard
scalar product $\langle \cdot,\cdot \rangle$ and by $\Ha^{\otimes
n}_q$ the $n$--fold tensor product equipped with the scalar
product $\langle \cdot,\cdot \rangle_q$.

The $q$--Fock space $\F$ is a complex Hilbert space defined by
$$ \F=\F(\Ha_\R)=\bigoplus_{n\geq 0} \Ha^{\otimes n}_q,$$
where the term $\Ha^{\otimes 0}_q$ should be understood as
one--dimensional space $\C \Omega$ for some unital vector
$\Omega$. By $\Ha \otimes \F$ we will mean the tensor product of
Hilbert spaces $\Ha$ and $\F$ equipped with the canonical scalar
product of the scalar product in $\Ha$ and the scalar product in
$\F$.

We consider the Hilbert space
$$ \F^{+}=\bigoplus_{n\geq 1} \Ha^{\otimes n}_q=\Omega^{\perp}  \subset \F$$ and the map
\begin{equation} \label{eq:zawieranie} j:\Ha
\otimes \F \rightarrow \F^{+}
\end{equation}
given by the trivial mapping of Hilbert spaces
$$ \phi\otimes (\psi_1\otimes \cdots \otimes
\psi_n)\mapsto \phi\otimes \psi_1\otimes \cdots \otimes \psi_n. $$
Please note that the map \eqref{eq:zawieranie} is not as trivial
as it might appear since the scalar products are different on the
domain and on the range.

\begin{proposition}
\label{prop:inclusion}
For every $-1<q<1$ there exist constants
$C_1$, $C_2$ such that for every choice of the real Hilbert space
$\Ha_\R$
$$\| j \|\leq C_1, \qquad \| j^{-1} \| \leq C_2. $$

Similar inequalities hold for the right trivial mapping
$\F\otimes\Ha\rightarrow\F^{+}$ given by $ (\psi_1\otimes \cdots
\otimes \psi_n)\otimes \phi\mapsto \psi_1\otimes \cdots \otimes
\psi_n\otimes\phi$.
\end{proposition}
\begin{proof}
Bo\.zejko and Speicher showed that
$$ P^{(n+1)} \leq \frac{1}{1-|q|} \big( 1 \otimes P^{(n)} \big) $$
and Bo\.zejko \cite{Bozejko1997} showed that there exists a
positive constant $\omega(q)$ such that for each $n$ we have
$$1 \otimes P^{(n)} \leq \frac{1}{\omega(q)}\ P^{(n+1)}. $$
\end{proof}

\subsection{Annihilation and creation operators}
For $\phi\in\Ha_\R$ we consider left and right creation operators
$\leftfreecre_\phi,\rightfreecre_\phi:\F \rightarrow \F$ defined
on elementary tensors by
\begin{equation}
\label{eq:leftqcre} \leftfreecre_\phi (\psi_1\otimes \cdots
\otimes \psi_n)=\phi\otimes\psi_1\otimes \cdots \otimes
\psi_n,\end{equation}
\begin{equation}
\label{eq:rightqcre}
\rightfreecre_\phi (\psi_1\otimes \cdots
\otimes \psi_n)=\psi_1\otimes \cdots \otimes \psi_n\otimes\phi.
\end{equation} We also consider
their adjoints: left and right
annihilation operators \begin{equation} \label{eq:leftqan}
\leftqan_\phi (\psi_1\otimes \cdots \otimes \psi_n)= \sum_{1\leq i
\leq n} q^{i-1} (\phi,\psi_i) \ \psi_1 \otimes \cdots \otimes
\widehat{\psi_i} \otimes \cdots \otimes \psi_n, \end{equation}
\begin{equation}
\label{eq:rightqan} \rightqan_\phi (\psi_1\otimes \cdots \otimes
\psi_n)= \sum_{1\leq i \leq n} q^{n-i} (\phi,\psi_i) \psi_1
\otimes \cdots \otimes \widehat{\psi_i} \otimes \cdots \otimes
\psi_n, \end{equation} where $\widehat{\psi_i}$ denotes an omitted
factor and $n\geq 1$. The case $n=0$ should be understood as
$\leftqan_\phi \Omega =\rightqan_\phi \Omega=0$. All creation and
annihilation operators are bounded.

\begin{remark}
Please note that in the definition of the annihilation operators
the natural extension of the bilinear scalar product
$(\cdot,\cdot)$ was used and not the sesquilinear scalar product
$\langle \cdot,\cdot \rangle$.
\end{remark}


\subsection{Von Neumann algebra}

Let $\Gamma_q$ denote the von Neumann algebra generated by the
family of selfadjoint operators $\leftfield_\phi:\F\rightarrow\F$
given by
$$\leftfield_\phi=\leftqan_\phi+\leftqcre_\phi.  $$
for $\phi\in\Ha_\R$.  Algebra $\Gamma_q$ is equipped with a
faithful tracial state $x\mapsto \langle\Omega, x\Omega\rangle$.
Vector $\Omega$ is separating and cyclic for $\Gamma_q$; let
$J:\F\rightarrow\F$ denote the canonical involution. Then $J
L_\phi J=R_\phi$, where
$$\rightfield_\phi=\rightqan_\phi+\rightqcre_\phi$$ for
$\phi\in\Ha_\R$.  Operators $(\leftfield_\phi)$ commute with
operators $(\rightfield_\psi)$ \cite{BozejkoXu}.



\section{The main result}
\label{sec:main} Let us fix some finite--dimensional real Hilbert
subspace $\Ha_\R'\subseteq\Ha_\R$ and let $e_1,\dots,e_d$ be an
orthonormal basis of $\Ha_\R'$.  We denote by $\Ha'$ the
complexification of $\Ha_\R'$.

In Section \ref{sec:inequalities} we will define a certain
operator $\cc$. It is enough to know now that
$|\cc|^2:\F\rightarrow\F$ is given by
\begin{equation} |\cc|^2=\cc\gwia \cc=\sum_i (L_{e_i}-R_{e_i})^2.
\label{eq:cekwadrat}
\end{equation}
In Section \ref{sec:inequalities} we will show the following.

\begin{proposition}
\label{prop:ugaugauga}
$\C\Omega$ belongs to the kernel of $| \cc
|$.

For every $-1<q<1$ there exists $d_0$ such that if $d\geq d_0$
then the restriction of $|\cc|$ to the space $\F^+$ is strictly
positive in the sense that
$$ |\cc| \geq \epsilon>0 $$
holds for some $\epsilon$.
\end{proposition}

The above result has an immediate consequence.

\begin{theorem}
\label{theo:main} Let $d_0$ be the constant from Proposition
\ref{prop:ugaugauga}. If $\dimm \Ha_\R \geq d_0$ (the case $\dimm
\Ha_\R=\infty$ is allowed) then $\Gamma_q$ is a $\text{II}_1$
factor which does not have the property $\Gamma$.
\end{theorem}
\begin{proof}
We choose finite--dimensional $\Ha_\R'\subseteq\Ha_\R$ in such a
way that $\dimm \Ha_\R'\geq d_0$.

Let $a\in\Gamma_q$ be central; it follows that $a$ commutes with
$| \cc|$ hence $ | \cc | a \Omega = a | \cc | \Omega = 0$ and $a
\Omega$ belongs to the kernel of $|\cc|$. Proposition
\ref{prop:ugaugauga} implies that $a \Omega \in \C \Omega$.
 Since $\Omega$ is separating it follows that
$a\in\C$ hence $\Gamma_q$ is a factor.

The following observation was pointed out to me by Dimitri
Shlyakhtenko. The operator $|\cc|$ belongs to the $C\gwia$-algebra
generated by $\Gamma_q$ and $J\Gamma_q J$. Proposition
\ref{prop:ugaugauga} shows that the kernel of $|\cc|$ is equal to
$\C \Omega$ and that the zero eigenvalue is separated. It follows
that the orthogonal projection $\F\rightarrow\C\Omega$ belongs to
$C\gwia(\Gamma_q,J\Gamma_q J)$. By the result of Connes
\cite{ConnesClassification} it follows that $\Gamma_q$ does not
have the property $\Gamma$.

\end{proof}


%


\section{Inequalities}
\label{sec:inequalities}

We define $\cc:\F \rightarrow \Ha \otimes \F$ by
$$\cc (\Phi)=\sum_i e_i\otimes \big(
L_{e_i}(\Phi)-R_{e_i}(\Phi) \big). $$ It is easy to check that
\eqref{eq:cekwadrat} is indeed fulfilled. We also define auxiliary
maps $\can,\ccre:\F \rightarrow \Ha \otimes \F$ given by
$$\can(\Phi)=\sum_i e_i\otimes \big(
\leftqan_{e_i}(\Phi)-\rightqan_{e_i}(\Phi) \big), $$
$$\ccre(\Phi)=\sum_i e_i\otimes \big( \leftqcre_{e_i}(\Phi)-\rightqcre_{e_i}(\Phi) \big). $$
Clearly $\cc=\can+\ccre$. Please note that the definitions of the
above operators do not depend on the choice of the orthonormal
basis $(e_i)$.

\begin{lemma}
\label{lem:malo} The norm of $\can$ fulfills
 $$\| \can \| \leq 2 C_1, $$
where $C_1$ is the constant from Proposition \ref{prop:inclusion}.
\end{lemma}
\begin{proof}
Consider the map $\can_l:\F^{+}\rightarrow \Ha  \otimes \F$ given
by
$$ \can_l \Phi =\sum_i e_i \otimes (l_{e_i} \Phi).$$
Its adjoint $\can_l \gwia:\Ha  \otimes \F \rightarrow \F^{+}$ is
given by
$$\can_l\gwia(\phi\otimes \Phi) = \sum_{i} (e_i,\phi)\
l^\ast_{e_i} \Phi =  j \big[ (\Pi_{\Ha'} \phi) \otimes \Phi
\big]$$ where $j$ denotes the trivial map \eqref{eq:zawieranie}
and $\Pi_{\Ha'}:\Ha\rightarrow\Ha'$ is the orthogonal projection.
Therefore $\|\can_l \|=\|\can_l\gwia\|\leq \|j\| \leq C_1$.

Similar inequality can be obtained for the right annihilator
$$ \can_r \Phi =\sum_i e_i \otimes (r_{e_i} \Phi);$$
it remains to notice that $\can=\can_l+\can_r$.
\end{proof}

%
%
%
%
%
%
%

\begin{lemma}
\label{lem:duzo} The restriction of $|\ccre|$ to the space $\F^+$
fulfils
$$ |\ccre| \geq   \frac{d-C_1 C_2}{C_2 \sqrt{d}}.$$
\end{lemma}
\begin{proof}
Consider the map $S:\F^{+} \rightarrow \F^{+}$:
$$S(\phi_1\otimes \cdots \otimes \phi_n)= \phi_2 \otimes \cdots
\otimes \phi_{n} \otimes (\Pi_{\Ha'} \phi_1) , $$ where
$\Pi_{\Ha'}:\Ha\rightarrow\Ha'$ denotes the orthogonal projection.
Proposition \ref{prop:inclusion} implies that $\|S\| \leq C_1
C_2$.



Consider the map $\tilde{f}:\Ha\otimes\Ha\otimes \F\rightarrow\F$
given by
\begin{multline*} \tilde{f}\big( \phi \otimes \psi_1
\otimes(\psi_2 \otimes \cdots \otimes \psi_n) \big)= \\
\Big\langle \sum_k e_k \otimes e_k, \phi \otimes \psi_1
\Big\rangle\ \psi_2 \otimes \cdots \otimes \psi_n;
\end{multline*}
its norm is equal to $\left\| \sum_k e_k \otimes e_k
\right\|=\sqrt{d}$. From Proposition \ref{prop:inclusion} it
follows that $f:\Ha\otimes \F^{+}\rightarrow\F$ given by
$$ f\big( \phi
\otimes (\psi_1\otimes \cdots \otimes \psi_n) \big)=
\Big\langle \sum_k e_k \otimes e_k, \phi \otimes \psi_1
\Big\rangle\ \psi_2 \otimes \cdots \otimes \psi_n
$$
fulfills $\|f\| \leq C_2 \sqrt{d}$ (because the only difference
between the maps $f$ and $\tilde{f}$ is the choice of the norm on
the domain).

It is easy to check that for $\Phi\in\F^{+}$
$$ f \ccre \Phi = (d - S) \Phi. $$

It follows
$$C_2 \sqrt{d}\ \|\ccre(\Phi)\| \geq \|f \ccre(\Phi) \|=\| (d-S)
\Phi \| \geq (d-C_1 C_2)\ \|\Phi\|.$$

%

\end{proof}

\begin{proof}[Proof of Proposition \ref{prop:ugaugauga}]
It is easy to check that $\Omega$ belongs to the kernel of
$|\cc|^2$.

Observe that if the restrictions of the operators to the space
$\F^+$ fulfill
$$ |\ccre|> \| \can\|  $$
then the restriction of $|\cc|=|\ccre+\can|$ to $\F^+$ is strictly
positive. Estimates from the above Lemmas finish the proof.
\end{proof}

\begin{remark}
For the case $q=0$ one can take $C_1=C_2=1$ and Theorem
\ref{theo:main} shows factoriality of $\Gamma_0$ for $\dimm
\Ha_\R\geq 6$ while we know (from the free group construction)
that it is valid also for $\dimm \Ha_\R\geq 2$ therefore our
result is far from being fully satisfactory.
\end{remark}

\section{Generalizations}

A careful reader might easily observe that in Sections
\ref{sec:main} and \ref{sec:inequalities} we used only very few
properties of Bo\.zejko--Speicher algebras. In particular, we did
not use the exact form of the symmetrizer \eqref{eq:symmetrizer}
and of the annihilation operators \eqref{eq:leftqan},
\eqref{eq:rightqan}, nor the existence of the tracial state
$a\mapsto\langle \Omega, a \Omega\rangle$.

It follows that our proof can be used in a more general context to
proof factoriality (both in type $\text{II}_1$ and type
$\text{III}$ case) of some von Neumann algebras arising from some
generalized Brownian motions for which an analogue of Proposition
\ref{prop:inclusion} holds.

It is easy to find a whole zoo of examples since already Bo\.zejko
and Speicher in their original article \cite{BozejkoSpeicher1994}
considered symmerizations arising from Yang--Baxter operators
which are much more general than \eqref{eq:symmetrizer}; Kr\'olak
\cite{Krolak2000} generalized further these results. Also Hiai
\cite{Hiai2002qdeformed} generalized
 results of Shlyakhtenko \cite{Shlyakhtenko1997quasifree} and
 constructed von Neumann algebras which are $q$--analogues of free
 Araki--Woods algebras. It was pointed out to me by M\u{a}d\u{a}lin
 Gu\c{t}\u{a} that the generalized Brownian motions considered in
 the work \cite{GutaMaassenGeneralized} might also provide appropriate
 examples.


\section{Acknowledgments}
I thank Marek Bo\.zejko for introducing me into the subject and
for many discussions. I thank Beno\^\i{}t Collins for teaching me
the right way of thinking about tensors. I also thank
M\u{a}d\u{a}lin Gu\c{t}\u{a}, Fumio Hiai, Dimitri Shlyakhtenko and
the Reviewer for many helpful remarks.

Research supported by State Committee for Scientific Research
(KBN) grant No.\ 2 P03A 007 23. The research was conducted in
Syddansk Universitet (Odense, Denmark) and Banach Center
(Warszawa, Poland) on a grant funded by European Post--Doctoral
Institute for Mathematical Sciences.

\def\cprime{$'$} \def\Dbar{\leavevmode\lower.6ex\hbox to 0pt{\hskip-.23ex
  \accent"16\hss}D} \def\cftil#1{\ifmmode\setbox7\hbox{$\accent"5E#1$}\else
  \setbox7\hbox{\accent"5E#1}\penalty 10000\relax\fi\raise 1\ht7
  \hbox{\lower1.15ex\hbox to 1\wd7{\hss\accent"7E\hss}}\penalty 10000
  \hskip-1\wd7\penalty 10000\box7}
  \def\cfudot#1{\ifmmode\setbox7\hbox{$\accent"5E#1$}\else
  \setbox7\hbox{\accent"5E#1}\penalty 10000\relax\fi\raise 1\ht7
  \hbox{\raise.1ex\hbox to 1\wd7{\hss.\hss}}\penalty 10000 \hskip-1\wd7\penalty
  10000\box7} \def\cprime{$'$} \def\cprime{$'$} \def\cprime{$'$}
  \def\cprime{$'$} \def\cprime{$'$} \def\cprime{$'$} \def\cprime{$'$}
  \def\cprime{$'$} \def\cprime{$'$} \def\cprime{$'$} \def\cprime{$'$}
  \def\cprime{$'$} \def\lfhook#1{\setbox0=\hbox{#1}{\ooalign{\hidewidth
  \lower1.5ex\hbox{'}\hidewidth\crcr\unhbox0}}}

\end{document}